\documentclass[12pt,reqno]{article}
\usepackage{amsmath, amsthm, graphicx,amsfonts, amssymb,color}
\usepackage{bm}
\usepackage{hyperref}

\usepackage{amsmath}
\numberwithin{equation}{section}
\def\d{\mathrm{d}}\def\e{\mathrm{e}}
\newtheorem{theorem}{Theorem}[section]
\newtheorem{lemma}[theorem]{Lemma}
\newtheorem{corollary}[theorem]{Corollary}

\newtheorem{remark}[theorem]{Remark}
\newtheorem{example}[theorem]{Example}

\def\<{\langle}\def\>{\rangle}
\def\prf{\noindent{\bf Proof.~~}}
\def\deprf{\hfill$\Box$\medskip}

\usepackage{mathrsfs}
\newcommand{\scr}[1]{\mathscr #1}

\newcommand{\set}[1]{\left\{#1\right\}}
\newcommand{\eps}{\varepsilon}

\def\d{\mathrm{d}}

\def\bg{\begin}

\def\de{\end{equation}}

\def\edar{\end{eqnarray}}

\def\lb{\label}

\def\l{\left}\def\r{\right}
\def\fr{\frac}
\def\alp{\alpha}
\def\bt{\beta}

\def\Gm{\Gamma}
\def\dlt{\delta}

\def\eps{\epsilon}

\def\kp{\kappa}
\def\lmd{\lambda}

\def\sgm{\sigma}

\def\ex{\exists}

\def\ift{\infty}

\def\rar{\rightarrow}

\def\q{\quad}

\def\[{\l[} \def\]{\r]}
\def\({\l(} \def\){\r)}

\def\bar{\overline}

\def\beqlb{\begin{eqnarray}}\def\eeqlb{\end{eqnarray}}
\def\beqnn{\begin{eqnarray*}}\def\eeqnn{\end{eqnarray*}}
\def\d{{\mbox{\rm d}}}

\title{{\bf Convergence rates in uniform ergodicity by hitting times and $L^2$-exponential convergence rates}\footnote{This work is supported in part by the National Natural Science Foundation of China (No.11771047) and the National Key Research and Development Program
of China (No. 2020YFA0712900).} }
\author{
{\bf Yong-Hua Mao and Tao Wang\footnote{Corresponding author: Tao wang
	(email: wang\_tao@mail.bnu.edu.cn; ORCID iD:	0000-0002-0727-5555)}}\\
\footnotesize{School of Mathematical Sciences, Beijing Normal University, }\\
\footnotesize{Laboratory of Mathematics and Complex Systems, Ministry of Education}\\
\footnotesize{Beijing 100875, China}\\
}
\date{ }

\date{}

\begin{document}

\maketitle


\bg{abstract}
Generally the  convergence rate in exponential ergodicity $\lmd$ is an upper bound for the convergence rate $\kp$
 in uniform ergodicity for a Markov process, that is $\lmd\geqslant\kp$. In this paper, we  prove that $\kp\geqslant \inf \set{\lmd,1/M_H}$, where $M_H$ is a uniform bound on the moment of the hitting time to a ``compact'' set $H$. In the case where $M_H$ can be made arbitrarily small for $H$ large enough, we obtain that $\lmd=\kp$. The general results are applied to Markov chains, diffusion processes and solutions to SDEs driven by symmetric stable processes.

\end{abstract}




\section{Introduction and general results}

Uniform ergodicity (or strong ergodicity) is an important topic in ergodic theory for Markov processes. In this paper, we are interested in the convergence rate in uniform ergodicity.

It is well known that the criterion for a Markov process to be uniformly ergodic, is to use the uniformly bounded moment of the first return time related to any petite set (or equivalently, a bounded Lyapunov function), especially for the Markov chains. See \cite{Ander91,DMPS18,MT95,MT09}.

To get the (exponential) convergence rates for discrete-time Markov chains, several types of classical methods  are  used,  such as minorization conditions (\cite{rosen95}), Foster-Lyapunov criteria (\cite{Bax05}) and  Dobrushin's ergodicity coefficients (\cite{sen88}) which can be used conceptually to continuous-time Markov processes, as in \cite[Chapter 6]{Ander91}. 


Coupling methods can be generally used to estimate the convergence rate via the moments of  the so-called coupling time (see \cite{cmf04}, \cite{cmf05}). This was done in \cite{myh02} for  the convergence rates in uniform ergodicity of Markov chain and diffusion process, then was improved by \cite{myh06}. However, to apply the coupling method, the stochastic monotone property is often technically needed to estimate the moment of the coupling time.

Historically, the study of the convergence rate in uniform ergodicity was much later than that in exponential ergodicity, although both the convergence rates are exponential. One reason for this may lie on the fact that $L^\infty$-norm for uniform ergodicity is less smooth than the  $L^2$-norm for exponentially ergodicity, especially for the reversible Markov processes. 
Even for the reversible Markov processes, no functional inequality can be adopted  directly for the convergence rates in uniform ergodicity.
For the reversible Markov processes, the spectral gap given by the classical Poincar\'{e} inequality is identical to the optimal convergence rate in exponential ergodicity (see \cite{cmf04}, \cite{cmf05} or \cite{wfy05}).

A ``mixed'' method appeared in \cite{myh10} where the moment of hitting time and spectral gap for reversible Markov chains are used to estimate the convergence rate in uniform ergodicity. 
The advantage of the ``mixed'' method is two-fold. 

On the one hand, in many cases, the 
uniform moment of hitting time can also afford the lower bound for the convergence rate in exponential ergodicity, so that we can get the explicit bounds  by using the moments of hitting times for many concrete models.

On the other hand, if it happens that the upper bound can be given by the  convergence rate in exponential ergodicity $\lambda$, then we find a phenomenon that the optimal convergence rate in uniform ergodicity $\kappa$ equals to $\lambda$ whenever the process is uniformly ergodic.
This is an interesting  phenomenon which was first proved in \cite{mz17} for the birth-death process. In general, if a reversible Markov semigroup $P_t$ is ultra-bounded, i.e. $\|P_t\|_{2\rightarrow\infty} <\infty$ for some $t>0$, then $\kappa=\lambda$, see \cite[Proposition 1.3]{myh06} for an argument. However, ultra-boundedness is a much stronger property to be satisfied.
As we will see soon, we actually find an extensive class of Markov processes, from Markov chains, diffusion processes to L\'{e}vy type processes, satisfying $\kappa=\lmd$.

Let $(X_t)_{t\geqslant 0}$ be a Markov process on state space $(E,\mathscr{B})$ with transition function $P_t(x,\cdot)$ which admits a stationary probability measure $\pi$.  

\bg{definition}\lb{def1}
The (exponential) {convergence rate in uniform ergodicity} is defined by
$$\kappa=
-\lim\limits_{t\rightarrow \infty}\frac{1}{t}\log \sup\limits_{x\in E}\|P_t(x,\cdot)-\pi\|_{\mathrm{Var}}.
$$
\end{definition}

If $X$ is uniformly ergodic, then $\sup_{x\in E}\|P_t(x,\cdot)-\pi\|_{\mathrm{Var}}\rar0$ as $t\rar\ift$. This convergence must be exponential, since by Markov property:
$$
\sup_{x\in E}\|P_{t+s}(x,\cdot)-\pi\|_{\mathrm{Var}}\leqslant\sup_{x\in E}\|P_t(x,\cdot)-\pi\|_{\mathrm{Var}}\times\sup_{x\in E}\|P_s(x,\cdot)-\pi\|_{\mathrm{Var}}.
$$
So $\ex C<\ift$ and $\eps>0$ such that 
$
\sup_{x\in E}\|P_t(x,\cdot)-\pi\|_{\mathrm{Var}}\leqslant C e^{-\eps t}.
$
Then $\kp$ is the maximal $\eps$ in the previous estimate.

Two basic concepts are related to our study on the convergence rate in uniform ergodicity. 

The first concept is the convergence rate in exponential ergodicity: 
there exist $\eps >0$ and a non-negative function $C(x)<\infty$ such that for any $x\in E$,
\begin{equation}\label{ee}
\|P_t(x,\cdot)-\pi(\cdot)\|_{\mathrm{Var}}\leqslant C(x) \mathrm{e}^{-\eps t}. 
\end{equation}
Denote by $\lambda$ the maximal $\eps$ in the above inequality, which is called the convergence rate in exponential ergodicity. 
Obviously, $\lambda\geqslant\kappa$.
A closed quantity to $\lmd$ is the $L^2$-exponential convergence rate $\lmd_1$:
$$
\lambda_1:=-\lim_{t\rightarrow \infty}\frac{1}{t}\log 
\|P_t-\pi\|_{L^2(\pi) \rar L^2(\pi)},
$$
where $ L^2(\pi)$ is the usual $L^2$-space with respective to $\pi$. 
For the reversible Markov processes, $\lambda_1$ is just the $L^2$-spectral gap:
$$
\lmd_1=\inf\{D(f,f):f\in \mathscr{D}, \pi(f)=0,\pi(f^2)=1\},
$$  
where $(D,\mathscr{D})$ is the Dirichlet form of $X$. 
In the reversible case,   denote by $p_{s}(x,y)$ the transition density with respect to $\pi$. If $p_{2s}(x,x) \in L_{\mathrm{loc}}^{1/2}(\pi)$ and the set of bounded
functions with compact support is dense in $L^2(\pi)$, then $\lambda=\lambda_1$ (cf.\cite[Theorem 8.13(4)]{cmf05}).  For the general Markov process, $\lambda$ and $\lambda_1$ may not be equal, but usually $\lmd\geqslant\lmd_1$ (see Corollary \ref{heat-ker} below).

The second concept related to $\kp$ is the uniform moment of hitting time:
$$
M_H:=\sup_{x\in E}\mathbb{E}_x \tau_H,
$$
where $\tau_H=\inf\{t\geqslant 0: X_t\in H\}$ is the hitting time to a  subset $H$. 
It is well-known that under some regular condition,  $X$ is uniformly ergodic if and only if  $M_H<\infty$ for some ``petite''  set $H$ (cf.\cite{Ander91, MT09} and reference therein).

In this paper, we will use exponential ergodicity convergence rate $\lmd$ or $\lmd_1$ and the moment $M_H$ to derive the convergence rate $\kp$ in uniform ergodicity. 
For this, {\bf unless otherwise stated,} we always make the following assumptions:

\begin{description}
\item[(A1)] The state space $(E,\mathscr{B})$ is a locally compact Polish space with metric $\rho$,  $X=(X_t)_{t\geqslant 0}$ is a progressive measurable  right continuous strong Markov process on a probability space $(\Omega, \mathscr{F}, (\mathscr{F}_t)_{t\geqslant0}, \mathbb{P})$, where $(\mathscr{F}_t)_{t\geqslant0}$ is the natural filtration;

\item[(A2)] $X$ is non-explosive, admits a stationary probability measure $\pi$.
\end{description}

Under Assumption (A1),  $\tau_H$ is a stopping time with respect to 
$(\mathscr{F}_t)_{t\geqslant0}$ and  $X_{\tau_H}\in H$ for non-empty closed set $H$.
Set
$$
\mathscr{H}=\{H\in \mathscr{B}: H \ \text{is a bounded closed set, such that}\ M_H=\sup_{x\in  E}\mathbb{E}_x \tau_H<\infty\}.
$$

Now we can claim our  general result, giving the relationship among $\kp$ and $\lmd, M_H$.

\begin{theorem}\label{main-thm}
Let $\lmd$ be the convergence rate in exponential ergodicity.
Assume that for any $\eps< \lmd$, 
(\ref{ee}) holds with $\sup_{x\in H}C(x)<\ift$ for some $H\in\mathscr H$. 
Then
\begin{equation}\label{low-bound}
\kappa\geqslant \min\left\{\lambda,\frac{1}{M_H}\right\}>0.
\end{equation}
Consequently,

(R1) if there exists $H\in \mathscr{H}$  with  $\sup_{x\in H}C(x)<\ift$ such that $ \lambda\leqslant 1/{M_H},$ then  $\kappa=\lambda;$ 

(R2) if there exists $H\in \mathscr{H}$  with  $\sup_{x\in H}C(x)<\ift$ such that $\lambda\geqslant  1/{M_H},$ then $\kappa\geqslant 1/{M_H}.$

\end{theorem}

To apply Theorem \ref{main-thm}, we need to prove the local boundedness of $C(x)$ on some $H\in \mathscr H$ in the exponential ergodicity (\ref{ee}).  For Markov chain,  we can consider $H$ as a single point and represent  $C(x)$ explicitly by stationary distribution $\pi$ (such as Example \ref{birth-death}). 
By using transition density $p_t(\cdot,\cdot)$ to represent $C(x)$, we can replace the exponential convergence rate $\lmd $ by the $L^2$-exponential convergence rate $\lmd_1$.

\begin{corollary}\label{heat-ker}
Assume that $P_t(x,\d y)=p_t(x,y)\pi(\d y), x,y\in E$.
If there is $s>0$ such that $\phi(x):=\|p_s(x,\cdot)\|_{L^2(\pi)}^2<\ift$, $\pi$-a.s.  (in the case of reversible processes, we have $\phi(x)=p_{2s}(x,x)$), then  $\lmd\geqslant\lmd_1$. If further $\sup_{x\in H}\phi(x)<\ift$ for some $H\in\mathscr H$, then 
$$
\kp\geqslant\min\set{\lmd_1,1/M_H}.
$$
\end{corollary}

\begin{remark}
{(1)}	According to \cite{myh06}, if $(X_t)_{t\geqslant 0}$ is uniformly ergodic, then 
$$\mathscr{G}:=\{B\in \mathscr{B}: B \ \text{is a bounded closed set, and}\ \pi(B)>0\}\subset \mathscr{H}.
$$
In fact, there may exist some set $H\in \mathscr{H}\setminus \mathscr{G}$. For example,  we can take $H$ a singleton, say $\set{0}$, for the one-dimensional  $\alpha$-stable process with $\alpha\in(1,2)$. Although $\pi(\set{0})=0$ as $\pi$ has density with respect to the Lebesgue measure,  $M_{\set{0}}$ can be represented explicitly for the  ergodic time-changed  $\alpha$-stable process (see Theorem \ref{scstable} below). 

(2) When   $ \lambda\leqslant 1/{M_H}$ (or $\lambda_1\leqslant 1/{M_H}$ for the reversible process), it is interesting to get that $\kappa=\lmd$ (resp. $\kappa=\lambda_1$), that is the convergence rates in uniform ergodicity and exponential ergodicity are identical. 
\end{remark}

To end this section, we would like to give two examples to illustrate that the situations (R1) and (R2) in Theorem \ref{main-thm}, respectively. 

\begin{example}\cite[Theorem 1.1]{mz17}\label{birth-death}
Let $(X_t)_{t\geqslant 0}$ be a birth-death process on $\mathbb Z_+$ with birth rates $b_i>0(i\geqslant0)$ and death rates $a_i>0(i\geqslant1)$. Assume the process has the 
$\infty$-entrance boundary in Feller's sense:
$$
\sum_{i=0}^\ift\mu_i\sum_{j=i}^\ift\fr1{\mu_j b_j}=\ift\q \text{and}\q S:=\sum_{i=0}^\ift\fr1{\mu_i b_i}\sum_{j=i+1}^\ift\mu_j<\ift,
$$
where $\mu_0=1, \mu_i=b_0\cdots b_{i-1}/a_1\cdots a_i (i\geqslant1)$.

The process is reversible with the stationary  distribution $(\pi_i)_{i\geqslant0}: \pi_i=\mu_i/\sum_{j=0}^\ift\mu_j.$ 
By using the coupling method and the stochastically monotone property, the estimate $\kp\geqslant1/(eS)$ was firstly given in \cite{myh02} and then was improved to $\kp\geqslant1/S$ in \cite{myh06}. 
Now this estimate is improved further in two ways  by  applying Theorem \ref{main-thm}. First,  we see that (\ref{ee}) holds with $\eps=\lmd_1$ and $C(x)=\sqrt{\pi_x^{-1}-1}$ for $x\in \mathbb{Z}_+$. By putting $H_n=\{0,1,\cdots, n\}$, we have $\lim_{n\rightarrow\infty}M_{H_n}=0$, so that $ \lambda_1\leqslant 1/ M_{H_n}$ for $n$ large enough.  Hence $\kappa=\lambda_1>0$. From \cite{cmf00'}, we have 
$$
\dlt^{-1}\geqslant\lmd_1\geqslant (4\dlt)^{-1},
$$
where $\dlt=\sup_{n\geqslant0}\sum_{i=n}^\ift\mu_i\sum_{i=0}^{n-1}\fr1{\mu_i b_i}$. Hence $\dlt^{-1}\geqslant\kappa\geqslant (4\dlt)^{-1}$. Moreover, the approximation procedure in \cite{cmf00'} can be applied to $\kp$.
Second, from Theorem \ref{mc1} in Section 4 below, we have
$$
\kp=\lmd_1\geqslant\sup_{i\geqslant0}\( \max\set{S_i,\bar S_i}\)^{-1}\geqslant 1/S,
$$
where $S_i=\sum_{k=i}^\ift\fr1{\mu_k b_k}\sum_{j=k+1}^\ift\mu_j $ and $\bar S_i=\sum_{k=0}^{i-1}\fr1{\mu_k b_k}\sum_{j=0}^k\mu_j$. 
\end{example}
We also remark that the uniform ergodicity can not imply the ultra-contraction. In \cite{WJ09}, the examples of the uniformly ergodic birth-death processes were given to exclude the hyper-contraction, let alone ultra-contraction.

The argument in Theorem \ref{main-thm} can be also applied to the discrete-time Markov chains, as shown in  the following example.

\begin{example}\cite[Theorem 1.4]{myh10}\label{dis-chain}
Let $(X_n)_{n\geqslant 0}$ be a reversible Markov chain on a discrete state space $E$, with nonnegative definite transition matrix $P$ and stationary distribution $\pi$. Let $H=\{0\}\subset E$ and $ M_0:=\sup_{i\in E}\mathbb{E}_i \tau_0<\infty$.
According to \cite[Lemmas 3.11-3.12]{ST89}, the spectral gap $\lambda_1\geqslant {1}/{ M_0},$ therefore $\kappa\geqslant {1}/{ M_0}$ by Theorem \ref{main-thm} (R2).
\end{example}


The paper is organized as follows. In Section \ref{relation}, we prove our main results, which establish the relation among $\kappa$, $\lambda_1$ and $M_H$ for general Markov processes and obtain a new estimate of lowed bound for  $\kappa$. 
In Section \ref{esti by spec} and \ref{esti by hit},  we study two typical situations which made $\kappa=\lambda$ and $\kappa\geqslant1/M_0$ respectively. The processes  include Feller process with non-negative jump, single death process, 
diffusion process on manifold, and  SDE driven by symmetric stable process.

\section{Proof of main results}\label{relation}

The following lemma is the start point of our method, which can be seen as a mixture of hitting time and exponential ergodicity.

\begin{lemma}\label{main-eq}
For $H\in \scr{H}$, let $F_{x,H}(t)=\mathbb{P}_x[\tau_H\leqslant t]$ be the 
distribution of $\tau_H$ and $f(x,t)= \|P_t(x,\cdot)-\pi\|_{\mathrm{Var}}$.
Then
\begin{equation*}
f(x,t)\leqslant  \mathbb{P}_x[\tau_H> t]+\int_{0}^{t}\sup_{y\in H}f(y,t-s)\d F_{x,H}(s), \q x\notin H.
\end{equation*}
\end{lemma}

\prf
For  $x\notin H$ and $A\in\mathscr{B}$, we have
\begin{equation}\label{total var}
\begin{split}
|P_t(x,A)-\pi(A)|=&\big|\mathbb{P}_x[X_t\in A, \tau_H> t]+\mathbb{P}_x[X_t\in A, \tau_H\leqslant t]-\pi(A)\big|\\
\leqslant&\big||\mathbb{P}_x[X_t\in A, \tau_H> t]-\pi(A)\mathbb{P}_x[\tau_H> t])\big|\\
&\q\q \q+\big|\mathbb{P}_x[X_t\in A, \tau_H\leqslant t]-\pi(A)\mathbb{P}_x[ \tau_H\leqslant t]\big|\\
\leqslant& \mathbb{P}_x[ \tau_H> t]+\big|\mathbb{P}_x[X_t\in A, \tau_H\leqslant t]-\pi(A)\mathbb{P}_x[ \tau_H\leqslant t]\big|,
\end{split}
\end{equation}
where in the second inequality we use the fact $|a-b|\leqslant c$ for  $0\leqslant a,b\leqslant c$.

Note that by strong Markov property, for any $A\in \mathscr{B}$,  on $\{\tau_H\leqslant t\}$,
\begin{equation}\label{str_Mkv}
\begin{split}
\mathbb{P}_x[X_t\in A |\mathscr{F}_{\tau_H}]&=\mathbb{P}_x[X_{t-\tau_H}\circ\theta_{\tau_H}\in A |\mathscr{F}_{\tau_H}] =\mathbb{P}_{X_{\tau_H}}[X_{t-\tau_H}\in A ],
\end{split}
\end{equation}
where $\theta_s$ is the usual shift operator such that $X_{s+t}=X_t\circ\theta_s$ for $s,t\geqslant0.$
Using the conditional expectation with respect to the stopping $\sgm$-algebra $\mathscr{F}_{\tau_H}$, it follows  from \eqref{str_Mkv} that
\begin{equation}
\begin{split}
\mathbb{P}_x[X_t\in A, \tau_H\leqslant t]&=\mathbb{E}_x\big[\mathbb{P}_x[X_t\in A, \tau_H\leqslant t|\mathscr{F}_{\tau_H}]\big]=\mathbb{E}_x\big[\mathbf{1}_{\{\tau_H\leqslant t\}}\mathbb{P}_x[X_t\in A |\mathscr{F}_{\tau_H}]\big]\\ &=\int_E\int_{0}^{t}P_{t-s}(y,A)\mathbb{P}_{x}(\tau_H\in \d s,X_{\tau_H}\in\d y).\\
\end{split}
\end{equation}
Since $X_{\tau_H}\in H$, we have
\begin{equation}\label{I_t}
\begin{split}
\big|\mathbb{P}_x[X_t\in A, \tau_H\leqslant t]-\pi(A)\mathbb{P}_x[ \tau_H\leqslant t]\big|		
&=\left|\int_E\int_{0}^{t}(P_{t-s}(y,A)-\pi(A))\mathbb{P}_{x}(\tau_H\in \d s,X_{\tau_H}\in\d y)\right|\\
&\leqslant\int_{0}^{t}\sup_{y\in H}|P_{t-s}(y,A)-\pi(A)|\d F_{x,H}(s).
\end{split}
\end{equation}
By combining \eqref{total var} and \eqref{I_t},  the desired result is obtained.

\deprf

Now, we use the above lemma to prove Theorem \ref{main-thm}.

\noindent\textbf{\textbf{Proof of Theorem \ref{main-thm}}}.

(a) Thanks to exponential ergodicity (\ref{ee}), the integral by parts gives for  $x\notin H$,
\begin{equation}
\begin{split}
\int_{0}^{t}\sup_{y\in H}f(y,t-s)\d F_{x,H}(s)&\leqslant C_H\int_{0}^{t}\mathrm{e}^{-\eps(t-s)}\mathrm{d}(-\mathbb{P}_{x}[\tau_H>s])\\
&=C_H\mathrm{e}^{-\eps t}\left(1-\mathrm{e}^{\eps t}\mathbb{P}_{x}[\tau_H>t]+\int_{0}^{t}\mathbb{P}_{x}[\tau_H>s]\epsilon\mathrm{e}^{\eps s}\d s\right)\\
&\leqslant C_H\mathrm{e}^{-\eps t}\left(1+\int_{0}^{t}\mathbb{P}_{x}[\tau_H>s]\eps\mathrm{e}^{\eps s}\d s\right),
\end{split}
\end{equation}
where $C_H:=\sup_{x\in H}C(x)$.

(b) By \cite[Lemma 3.7]{cz95},
\begin{equation*}
\sup_{x\in E}\mathbb{E}_x[\tau_H^n]\leqslant n!M_H^n, \q \text{for}\  n=0,1,2,\cdots
\end{equation*}
so that 
\begin{equation}\label{exp mom}
\mathbb{E}_x[\e^{\beta\tau_H}]=\sum_{n=0}^\ift\fr{\beta^n\mathbb{E}_x[\tau_H^n]}{n!}\leqslant \frac{1}{1-\beta M_H},\q \text{ for $0<\beta<1/M_H$.}
\end{equation}
Thus
\begin{equation}\label{tail_prob}
\mathbb{P}_x[\tau_H>t]\leqslant \mathbb{E}_x[\e^{\beta\tau_H}]\e^{-\beta t}\leqslant \frac{1}{1-\beta M_H}\e^{-\beta t}.
\end{equation}

By (a), we have for $\eps\neq\beta$,
\begin{equation}
\begin{split}
\int_{0}^{t}\mathbb{P}_{x}[\tau_H>s](\eps\mathrm{e}^{\eps s})\d s  &\leqslant \frac{1}{1-\beta M_H}\int_0^t \mathrm{e}^{-\beta s} \mathrm{e}^{\eps s}\d s= \frac{\mathrm{e}^{(\eps-\beta) t}-1}{(\eps-\beta)(1-\beta M_H)},\\
\end{split}
\end{equation}
where in the case of $\beta=\eps$, the last term is understood as the limit of  $\beta\rightarrow\eps$.

(c) From (a) and (b), it follows that for $x\notin H$,
\begin{equation}
\begin{split}
f(x,t)&\leqslant \frac{2\mathrm{e}^{-\beta t}}{1-\beta M_H}+C_H\e^{-\eps t}
\left(1+\frac{\mathrm{e}^{(\eps-\beta) t}-1}{(\eps-\beta)(1-\beta M_H)}\right),\\
\end{split}
\end{equation}
while obviously for $x\in H$,
$$ f(x,t) \leqslant C_H\e^{-\eps t}.
$$
Therefore we have
$$\kappa= -\lim\limits_{t\rightarrow \infty}\frac{1}{t}\log \sup\limits_{x\in E}f(x,t)\geqslant\min\left\{\eps,\beta\right\}
$$ 
for any $\beta<1/M_B$ and $\eps<\lmd$, so that 
$$\kappa\geqslant \min\left\{\lambda,\frac{1}{ M_H}\right\}. $$
\deprf

\noindent\textbf{\textbf{Proof of Corollary \ref{heat-ker}}}.

(a) In the reversible case, we have $p_t(x,y)=p_t(y,x)$, $\pi\times\pi$ a.s. $(x,y)$,  hence  $\phi(x)=\|p_s(x,\cdot)\|_{L^2(\pi)}^2=p_{2s}(x,x)$. Since  $\lambda_{1}$ is equal to the spectral gap, $\lambda\geqslant \lambda_1$ by  \cite[Theorem 8.8]{cmf05}. 

(b) Denote by $\pi:f\mapsto\pi(f):=\int_E f\d \pi$. For the general case, by definition of $\lmd_1$, for any $\eps<\lmd_1$, there is $C_1<\ift$ such that 
$$
||P_t-\pi||_{L^2(\pi)\rar L^2(\pi)}\leqslant C_1 e^{-\eps t}, \q  t\geqslant0.
$$
Let $P_t^*$ is the dual semigroup of $P_t$ with respect to $\pi$.
Then
$$
||P_t^*-\pi||_{L^2(\pi)\rar L^2(\pi)}=||P_t-\pi||_{L^2(\pi)\rar L^2(\pi)}\leqslant C_1 e^{-\eps t}, \q  t\geqslant0.
$$ 

(c)  For $t\geqslant s$,  we have 
\begin{equation}
\begin{split}
P_tf(x)& =P_sP_{t-s}f(x)=\int_{E}p_s(x,y)P_{t-s}f(y)\pi(\d y)\\
&=\int_{E}[P_{t-s}^*(p_s(x,\cdot))(y)]f(y)\pi(\d y).
\end{split}
\end{equation}
So  by Cauchy-Schwartz inequality and (b),
\begin{equation}
\begin{split}
\|P_t(x,\cdot)-\pi\|_{\mathrm{Var}}&=\sup_{|f|\leqslant 1}\big|P_tf(x)-\pi(f)\big| \\
&=\sup_{|f|\leqslant 1}\left|\int_{E}[P_{t-s}^*(p_s(x,\cdot))(y)-1]f(y)\pi(\d y)\right|\\
&\leqslant \|P_{t-s}^*\big(p_s(x,\cdot)-1\big)\|_{L^2(\pi)}\\
&\leqslant \|p_s(x,\cdot)-1\|_{L^2(\pi)}C_1 \e^{-\eps(t-s)}.
\end{split}
\end{equation}
Since $\|P_t(x,\cdot)-\pi\|_{\mathrm{Var}}\leqslant2$, (\ref{ee}) holds for all $t\geqslant0$ by choosing $C(x)=\max\{2,\|p_s(x,\cdot)-1\|_{L^2(\pi)}C_1 \e^{\eps s}\}$.
Thus  $\lmd\geqslant\eps$  for any  $\eps<\lmd_1$, so that $\lmd\geqslant\lmd_1$. Then the desired result follows from Theorem \ref{main-thm}.

\deprf

Next, Sections \ref{esti by spec} and \ref{esti by hit} will discuss two classes of models which respectively satisfy the assumptions (R1) and (R2) in Theorem \ref{main-thm}.

\section{Estimate of $\kappa$ by  $\lmd$ or $\lmd_1$}\label{esti by spec}

In this section, we will seek the situation that the moment of hitting time can be used as the upper bound for $\lambda_1$ or $\lambda$.  More technically, we give the conditions for which $M_H$ goes to zero when $H$ becomes bigger and bigger. 
More generally, if there exists an $H\in\scr{H}$, such that 
$$\lambda\leqslant\frac{1}{M_H},\ \text{or,}\ \lambda_1\leqslant\frac{1}{M_H},\q \text{respectively,}
$$
then it holds that 
$$\kappa=\lambda \ \text{or}\ \kappa=\lambda_1 \q \text{in the reversible case respectively}.
$$
As have done in Example \ref{birth-death} for the birth-death process,   we  will do this by seeking a sequence $\{H_n\}\subset\scr{H}$ such that $H_n\uparrow E$ and 
\begin{equation}\label{instan-hit}
\lim_{n\rightarrow\infty}\sup_{x\in E}\mathbb{E}_x\tau_{H_n}=0,
\end{equation}
So there exists $H_n\in\scr{H}$ such that $\lambda_1\leqslant1/M_{H_n}$ or $\lambda\leqslant1/M_{H_n}$.

In the following subsections,   to study this situation, we present
a class of models including Markov processes with $\infty$ instantaneous entrance boundary, Markov chains, diffusion processes and   SDEs driven by symmetric stable processes.

\subsection{Feller processes with non-negative jumps}

Let $E=[0,\infty)$ and   $X$ be a  non-explosive Feller process with non-negative jump on $E$. 
We say $\infty$ is an {\bf instantaneous entrance boundary}, if  for any $t> 0$,
\begin{equation}\label{instan-entr1}
\lim_{b\rightarrow\infty}\limsup_{x\rightarrow\infty}\mathbb{P}_x(\tau_{[0,b]} > t) = 0.
\end{equation}
Cf. \cite{FLZ20}. It is proved in 	\cite[Lemma 1.2]{FLZ20} that for this process, \eqref{instan-entr1} is equivalent to 
$$\lim_{b\rightarrow\infty}{\lim_{x\rightarrow\infty}}\mathbb{E}_x\tau_{[0,b]}=0.
$$
By \cite[Proof of Lemma 1.2]{FLZ20}, for any $x>x'>b>0$,
$$\mathbb{E}_{x}\tau_{[0,b]}=\mathbb{E}_{x}\tau_{[0,x^{\prime}]}+\mathbb{E}_{x^{\prime}}\tau_{[0,b]}\geqslant\mathbb{E}_{x^{\prime}}\tau_{[0,b]},$$
i.e. $\mathbb{E}_{x}\tau_{[0,b]}$ is non-decreasing for $x>0$. Thus ${\lim_{x\rightarrow\infty}}\mathbb{E}_x\tau_{[0,b]}={\sup_{x}}\mathbb{E}_x\tau_{[0,b]}$, so \eqref{instan-hit} holds.
This ensures that 
$\kappa=\lambda$  by Theorem \ref{main-thm}.

\subsection{Single death processes}

As a counterpart of Markov process on $[0,\infty)$ with no negative jump, we consider the so-called single death process (or downwardly skip free process) on $\mathbb{Z}_+$. 

The $Q$-matrix $Q=(q_{ij})_{i,j\in\mathbb{Z}_+}$ is  called a single death $Q$-matrix, if $q_{i,i-1}>0 $ for all $i\geqslant1$, and $q_{i,i-j}=0$ for $i\geqslant j\geqslant2.$ Assume that $Q$ is regular, i.e.
$$q_{i}:=-q_{i i}=\sum_{j \neq i} q_{i j}<+\infty, \quad i \in \mathbb{Z}_{+},$$
and irreducible.
Let
$$q_n^{(k)}=\sum_{j=k}^{\infty}q_{nj},\ k>n\geqslant0,
$$
and define inductively 
$$  G_n^{(n)}=1, \ \ G_n^{(i)}=\frac{1}{q_{n,n-1}}\sum_{k=n+1}^{i}q_n^{(k)}G_k^{(i)}, \ \ 1\leqslant n<i.
$$

It is proved in \cite[Lemma 2.7]{zyh18} that the single death process is uniformly ergodic if and only if
$$
S:=\sum_{k=1}^{\infty}\sum_{l=k}^{\infty}\frac{G_k^{(l)}}{q_{l,l-1}}<\infty.
$$
Furthermore, for $i> n,$
$$
\mathbb{E}_i\tau_{n}\leqslant\sum_{k=n+1}^{i}\sum_{l=k}^{\infty}\frac{G_k^{(l)}}{q_{l,l-1}},
$$
where $\tau_n:=\inf\{t\geqslant0: X_t=n\}$. 
By choosing 
$$H_n=\{0,1,\cdots,n\},
$$
we have $\mathbb{E}_i\tau_{H_n}=\mathbb{E}_i\tau_{n}$ for $i>n$ by skip free property,  so that
$$
\sup_{i> n }\mathbb{E}_i\tau_{H_n}\leqslant \sum_{k=n+1}^{\ift}\sum_{l=k}^{\infty}\frac{G_k^{(l)}}{q_{l,l-1}}\rar 0, \q\text{as}\q n\rar\ift,
$$
provided $S<\ift$. Then 
$\kappa=\lambda>0$  by applying Theorem \ref{main-thm}.


\subsection{Diffusions processes}

First we consider the one-dimensional diffusion process which is both stochastically monotone Markov process and Feller process with non-negative jump.
\begin{corollary}[Diffusions on half-line]\label{1d diff}
Let $L=a(x)\frac{\mathrm{d}^2}{\mathrm{d}x^2}+b(x)\frac{\mathrm{d}}{\mathrm{d}x}$ be a diffusion operator on $\mathbb{R}_+$ with $a(x)>0$, and $a,b$ be continuous. Define $c(x)=\int_{1}^{x}\frac{b(y)}{a(y)}\mathrm{d}y,$ and $\pi(\d z)=a(z)^{-1}\mathrm{e}^{c(z)}\mathrm{d}z.$ Denote by $X$ the diffusion process on $[0,\infty)$ with generator $L$ and reflecting boundary at 0.
Assume that  $X$ has $\infty$-entrance boundary:
\begin{equation*}
\int_{0}^{\infty}\mathrm{e}^{-c(y)}\left(\int_{0}^{y}\frac{\mathrm{e}^{c(z)}}{a(z)}\mathrm{d}z\right)\mathrm{d}y=\infty,
\q 		\int_{0}^{\infty}\mathrm{e}^{-c(y)}\left(\int_{y}^{\infty}\frac{\mathrm{e}^{c(z)}}{a(z)}\mathrm{d}z\right)\mathrm{d}y<\infty.
\end{equation*}
Then $\kappa=\lambda_1=\inf\left\{\pi(a(f')^2):\pi(f)=0,\pi(f^2)=1\right\}. $
\end{corollary}
\prf
By \cite[Section 4.11]{IM65},  the heat kernel $p_t(x,x)$ can be chosen to be continuous in $x\in \mathbb{R}_+$, so by Corollary \ref{heat-ker}, $C(x)$ is locally bounded.
Note that (cf. \cite[Proof of Theorem 2.1]{myh02})
$$ M_r:=\sup\limits_{x>r}\mathbb{E}_x \tau_{[0,r]}=\int_{r}^{\infty}\mathrm{e}^{-c(y)}\left(\int_{y}^{\infty}\frac{\mathrm{e}^{c(z)}}{a(z)}\mathrm{d}z\right)\mathrm{d}y<\infty.
$$
Then $\lim\limits_{r\rightarrow\infty} M_r=0$, 
so by Theorem \ref{main-thm}  we have $\kappa=\lambda_1. $
\deprf

The above result provides a way by using the spectral gap $\lmd_1$ to estimate $\kp$ for the one-dimensional diffusion process with entrance boundary. For examples, the following estimate in {\cite{cmf99}} can be served as the estimate for $\kp$:
$$
\dlt^{-1}\leqslant\kp=\lmd_1\leqslant (4\dlt)^{-1},
$$
where 
$$
\dlt=\sup_{x>0}\int_{0}^{x}\mathrm{e}^{-c(y)}\mathrm{d}y\int_{x}^{\infty}\frac{\mathrm{e}^{c(z)}}{a(z)}\mathrm{d}z<\infty.
$$
This estimate improves the estimate $\kp\geqslant 1/M_0$ in \cite{myh06} by using the coupling method. Moreover, in {\cite{cmf99}}, the approximation procedure of $\lambda_{1}$ now can also be applied to estimate $\kp$.


Next we turn to diffusion processes on manifolds.
Let $M$ be a connected, complete Riemannian manifold with empty boundary or convex boundary, and $(X_t)_{t\geqslant 0}$ be a non-explosive diffusion process on $M$ generated by $L=\Delta+Z$ with invariant probability measure $\pi$ (cf. see \cite[Theorem 3.1]{brw01} for the sufficient condition for the existence of invariant measure), where $\Delta$ is the Laplacian and $Z$ denotes both the $C^1$ vector field on $M$ and the corresponding derivative operator. Assume that the curvature condition is satisfied, i.e. there exists a constant $K$ such that $\mathrm{Ric}(Y,Y)-\langle\nabla_YZ,Y\rangle\geqslant -K\|Y\|^2$. 

Under these assumptions, the dimensional free Harnack inequality holds (see \cite[Theorem 2.3.3]{wfy13}), thus by \cite[Corollary 3.1(2)]{wy11}, there exists density $p_t(x,y)$ with respect to $\pi$. 

Let $\rho$ be the Riemannian metric. Fix a point $o\in M$, 
set $\rho(x)=\rho(o,x)$ and $D=\sup_{x\in M} \rho(x)$. Assume that $\mathrm{cut}(o)=\emptyset$. 
Fix $r_0>0,$ let
\begin{equation*}
\overline{C}(r)=\int_{r_0}^{r}\overline{\beta}(s)\d s,\ \ \ \overline{\beta}(r)\geqslant\sup_{\rho(x)=r}L\rho(x) \q \text{for} \ r>r_0,
\end{equation*}
and
\begin{equation}\label{rho-entra}
\overline{\delta}(\rho):=\int_ {0}^{D}\mathrm{e}^{-\overline{C}(y)}\left(\int_ {y}^{D}\mathrm{e}^{\overline{C}(z)}dz\right)dy.
\end{equation}

\begin{theorem}\label{non-sym}
If $(X_t)_{t\geqslant 0}$ is non-explosive and 
$\overline{\delta}(\rho)<\infty,$
then the convergence rate
$\kappa\geqslant\lambda_1.$ 
Specially in the reversible case,  i.e. $Z=\nabla V\cdot \nabla$ for some $V\in C^2(M)$, we have $$\kappa=\lambda_1=\inf\{\pi(|\nabla f|^2):\pi(f)=0,\pi(f^2)=1\},$$
where $\pi(\d x)=\e^{V(x)}\d x/\int_M\e^{V(x)}\d x$.
\end{theorem}

Before starting the proof of Theorem \ref{non-sym}, we need the following lemma whose proof is similar to  that of \cite[Theorem 2.4.4]{wfy13}.

\begin{lemma}\label{non-sym density}
Let  $p_t(x,y)$ is the transition density.  Then for any $s, r>0,$ and $x\in M$,
\begin{equation}\label{heat-est}
\|p_s(x,\cdot)\|_{L^2(\pi)}^2\leqslant \frac{1}{\pi(B(x,r))}\e^{U_s(r)}, 
\end{equation}
where $B(x,r)=\set{y\in M: \rho(x,y)\leqslant r}$ and $U_s(r)={Kr^2}/(\e^{2Ks}-1)$.
\end{lemma}
\prf
Let $p=2$ in dimension-free Harnack inequality
(see \cite[Theorem 2.3.3]{wfy13}),
we have  for any positive bounded function $f$, 
$$(P_sf)^2(x)\leqslant P_sf^2(y)\e^{U_s(\rho(x,y))}.$$
Hence
\begin{equation*}
\pi(f^2)=\pi P_sf^2\geqslant (P_sf)^2(x)\int_M \e^{-U_s(\rho(x,y))}\pi(\d y)\geqslant(P_sf)^2(x) \e^{-U_s(r)}\pi(B(x,r)).
\end{equation*}
By choosing $f(y)=n\wedge p_s(x,y)$,   we obtain that
\begin{equation*}
\left(\int_M(n\wedge p_s(x,y)) p_s(x,y)\pi(\d y)\right)^2\leqslant \frac{1}{\pi(B(x,r))}\e^{U_s(r)}\pi((n\wedge p_s(x,\cdot))^2).
\end{equation*}
Since 
$$
\int_M(n\wedge p_s(x,y)) p_s(x,y)\pi(\d y)\geqslant \pi((n\wedge p_s(x,y))^2),
$$ 
we have
\begin{equation*}
\pi((n\wedge p_s(x,y))^2)\leqslant  \frac{1}{\pi(B(x,r))}\e^{U_s(r)},
\end{equation*}
By letting $n\rightarrow\infty$, we get  \eqref{heat-est}. \deprf

\noindent\textbf{\textbf{Proof of Theorem \ref{non-sym}}}

Let
$$
u_p(r)=\int_ {p}^{r}\mathrm{e}^{-\overline{C}(y)}\left(\int_ {y}^{D}\mathrm{e}^{\overline{C}(z)}dz\right)dy
$$
and $\overline{\delta}_p(\rho)=\lim\limits_{r\rightarrow D}u_p(r)$.
Then $\overline{\delta}_p(\rho)<\infty$ and $u_p$ satisfies that $u_p''(r)+\overline{\beta}(r)u_p'(r)=-1.$
Hence for $x\in M$ with $\rho(x)=r$,
\begin{equation}\label{lyapunov1}
L[u_p\circ\rho](x)=u_p''[\rho(x)]+L\rho(x)u_p'[\rho(x)]\leqslant -1.
\end{equation}
Taking $f_p(x)=u_p\circ\rho(x)$, and ${B_p}=\{x\in M: \rho(x)\leqslant p\}$, by the well-posedness of martingale problem, we have
\begin{equation}
\mathbb{E}_x[f_p(X_{ t\wedge\tau_{B_p}})]-f_p(x)=\mathbb{E}_x\left[\int_{0}^{t\wedge\tau_{B_p}}Lf_p(X_s)\d s\right]\leqslant -\mathbb{E}_x[ t\wedge\tau_{B_p}],
\end{equation}
Note that $\mathbb{E}_x[f_p(X_{\tau_{B_p}})]=0$ and $\sup\limits_{x\notin {B_p}}f_p(x)=\overline{\delta}_p(\rho)$. By letting $t\rightarrow\infty$,  we have that
$$ 
M_p:=\sup\limits_{x\notin {B_p}}\mathbb{E}_x[\tau_{B_p}]\leqslant \overline{\delta}_p(\rho)<\infty.
$$
Hence $\lim\limits_{p\rightarrow D} M_p=0$.  According to Lemma \ref{non-sym density}, 
$\|p_s(x,\cdot)\|_{L^2(\pi)}$ is locally bounded, consequently  $\kappa\geqslant\lambda_1$ by Corollary \ref{heat-ker}. Specially, if $Z=\nabla V$ for some $V\in C^2(M)$, then the process is reversible with respect to $\pi$, so $\kappa=\lambda_1=\inf\{\pi(|\nabla f|^2):\pi(f)=0,\pi(f^2)=1\}.$
\deprf

Theorem \ref{non-sym} can improve the estimates in \cite{myh06} for the diffusion processes on $M$ by using coupling method.
Here we show an example:
\begin{example}\cite[Example 1.9]{cmf95}
Let $(X_t)_{t\geqslant 0}$ be a diffusion process on $\mathbb{R}^n$ with generator $L=\Delta+\nabla V\cdot \nabla$, $V=-|x|^4$.
By \cite[Example 3.6]{myh06}, we get a lower bound of $\kappa:$
\begin{equation*}
\begin{split}
\kappa\geqslant \frac{1}{\delta(M)}&=: \left(\frac{1}{4}\int_{0}^{\infty}\mathrm{e}^{(r/2)^4}\mathrm{d}r\int_{r}^{\infty}\mathrm{e}^{-(s/2)^4}\mathrm{d}s\right)^{-1}.
\end{split}
\end{equation*}
By \cite[Example 1.9]{cmf95}, we have
\begin{equation}
\delta(M)\leqslant \Gamma(5/4)+\frac{1}{8}.
\end{equation}
then $\kappa\geqslant (\delta(M))^{-1}\approx 0.9695.$

But on the other hand, it is obvious that $(X_t)_{t\geqslant 0}$ satisfies the condition of Theorem \ref{non-sym}, hence $\kappa=\lambda_1$.
In \cite[Example 4.11]{cmf97}, apply $I$-operator to $f(x)=\log (1+x)$ to derive
\begin{equation*}\kappa=\lambda_1 \gtrapprox 2.4395.\end{equation*}
\end{example}
\subsection{\bf SDEs driven by symmetric stable processes}
Let $(Z_t)_{t\geqslant 0}$ be a $d$-dimensional symmetric $\alpha$-stable process with generator $-(-\Delta)^{\alpha/2},$ which has the following expression:
$$-(-\Delta)^{\alpha/2}f(x):=\int_{\mathbb{R}^d\setminus\{0\}}\left(f(x+z)-f(x)-\nabla f(x)\cdot z\mathbf{1}_{\{|z|\leqslant 1\}}\right)\frac{C_{d,\alpha}}{|z|^{d+\alpha}} \d z,$$
where 
\begin{equation}\label{Cd,alpha}
C_{d,\alpha}=\frac{\alpha2^{\alpha-1}\Gamma((d+\alpha)/2)}{\pi^{d/2}\Gamma(1-\alpha/2)}
\end{equation}
is the normalizing constant.

Consider the following stochastic differential equation (SDE) driven by $\alpha$-stable process on $\mathbb{R}^d$:
$$\d X_t = \d Z_t+  b(X_t) \d t , \  X_0 = x,$$
where  $b: \mathbb{R}^d \rightarrow  \mathbb{R}^d$ is a continuous function satisfying  there exists a constant $\eta > 0$ such that for all $x, y\in\mathbb{R}^d,$
\begin{equation}\label{dissipation}
\<b(x)-b(y), x - y\> \leqslant \eta|x -y|^2.
\end{equation}

Under the condition \eqref{dissipation}, the SDE has the unique strong solution $(X_t)_{t\geqslant 0}$ which is 
strong Feller and Lebesgue irreducible (see e.g. \cite{wj13}).
By \cite{wj13},  the extended generator  $(L,D_w(L))$ is given as follows:
\begin{equation}\label{D_w}
D_w(L):=\left\{ f \in C^2(\mathbb{R}^d): \int_{\{|z|>1\}}[f(x+z)-f(x)]\frac{1}{|z|^{d+\alpha}}\d z<\infty, \ \text{for} \ x\in \mathbb{R}^d \right\},
\end{equation}
and for any $ \ f\in D_w(L),$ 
$$Lf(x)=\int_{\mathbb{R}^d\setminus\{0\}}\left(f(x+z)-f(x)-\nabla f(x)\cdot z\mathbf{1}_{\{|z|\leqslant 1\}}\right)\frac{C_{d,\alpha}}{|z|^{d+\alpha}} \d z+\<b(x),\nabla f(x)\>.$$

\begin{theorem}\label{stable}
Assume that $\alpha\in (1,2)$. If 
for some $\delta>1$, the following drift condition holds:
$$\<x,b(x)\>\leqslant -K|x|^{1+\delta},$$
then the process is uniformly ergodic, and  the convergence rate $\kappa\geqslant\lambda_1>0.$
\end{theorem}
\prf
Define
\begin{equation}\label{g}
g(r)=\inf\limits_{|x|\geqslant r}\left\{-\frac{\<x,b(x)\>}{|x|^2}\vee 0\right\}\geqslant Kr^{\delta-1},
\end{equation}
and   $$\tilde{g}(r)=\frac{1}{r}\int_{0}^{r}g(s)\d s\geqslant \frac{K}{\delta }r^{\delta-1}. $$
Obviously, $g(r)$ is a non-decreasing function, so that $\tilde{g}(r)\leqslant g(r)$, and \begin{equation}\label{integ-cond}
\delta_r:=\int_{r}^{\infty}\frac{1}{s\tilde{g}(s)} \d s\leqslant\frac{\delta}{K(\delta-1)}r^{1-\delta}<\infty.
\end{equation} 
For any  $r>0$, we  take nonnegative  function $u_r(x)\in C^2(\mathbb{R}^d)$ such that  for $|x|> r,$ $u_r(x)=r^{-1}+\int_{r}^{|x|}(s\tilde{g}(s))^{-1} \d s$ and for $|x|\leqslant r$, $u_r(x)\leqslant r^{-1}$. Then  $ u_r\leqslant \delta_r+r^{-1}=:\eta_r<\infty$, hence $u_r$ is bounded. 
Therefore,
\begin{equation}\label{>1}
\int_{\{|z|>1\}}[u_r(x+z)-u_r(x)]\frac{1}{|z|^{d+\alpha}}\d z\leqslant 2\eta_r \Gm_d\int_{1}^{\infty}\frac{1}{r^{1+\alpha}}\d r=\frac{2\eta_r\Gm_d}{\alpha}<\infty,
\end{equation}
where $\Gm_d=2\pi^{d/2}/\Gamma(d/2)$ is the volume of the sphere in $\mathbb{R}^d$, so  by \eqref{D_w}, $u_r(x)\in D_w(L).$
A direct computation shows that 
\begin{equation}\label{drift}
\<b(x),\nabla u_r(x)\>=\frac{\<x,b(x)\>}{|x|^2}\frac{1}{\tilde{g}(|x|)}\leqslant -\frac{g(|x|)}{\tilde{g}(|x|)}\leqslant-1,
\end{equation}
and for $|z|\leqslant 1$,
\begin{equation*}
\begin{split}
u_r(x+z)-u_r(x)-\<z,\nabla u_r(x)\>
&=\frac{1}{2}\<z,D^2u_r(\xi)z\>\\
&=\frac{1}{2}\left(\frac{|z|^2}{|\xi|^2\tilde{g}(|\xi|)}-\frac{\<z,\xi\>^2}{|\xi|^4\tilde{g}(|\xi|)}-\frac{\<z,\xi\>^2g(|\xi|)}{|\xi|^6(\tilde{g}(|\xi|))^2}\right)\\
\end{split}
\end{equation*}
where $\xi=x+\theta z$, $\theta\in(0,1).$
Note that when $|x|>1$ and $|z|\leqslant 1$,
$$|\xi|\geqslant|x|-\theta|z|\geqslant|x|-1.$$
Thus $$u_r(x+z)-u_r(x)-\<z,\nabla u_r(x)\>\leqslant \frac{1}{2(|x|-1)^2\tilde{g}(|x|-1)}|z|^2.$$
Therefore,
\begin{equation}\label{<1}
\begin{split}
\int_{\{|z|\leqslant1\}}[u_r(x+z)-u_r(x)-\<z,\nabla u_r(x)\>]\frac{C_{d,\alpha}\d z}{|z|^{d+\alpha}}&\leqslant \frac{1}{2(|x|-1)^2\tilde{g}(|x|-1)}\int_{\{|z|\leqslant1\}}\frac{|z|^2C_{d,\alpha}\d z}{|z|^{d+\alpha}}\\
&= \frac{ C_{d,\alpha}\Gm_d}{2(2-\alpha)(|x|-1)^2\tilde{g}(|x|-1)},
\end{split}
\end{equation}
where $C_{d,\alpha}$ is defined in \eqref{Cd,alpha}.
Combining \eqref{>1}, \eqref{drift} and \eqref{<1}, we get that  for $|x|\geqslant r$ and $r>1$,
$$Lu_r(x)\leqslant-1+\frac{2\eta_r\Gm_d}{\alpha}+\frac{ C_{d,\alpha}\Gm_d}{2(2-\alpha)(r-1)^2\tilde{g}(r-1)},
$$
so that  $Lu_r(x)\leqslant-\frac{1}{2}$ for $|x|\geqslant r$ with $r$ large enough.   By the definition of $D_w(L),$
\begin{equation}
\mathbb{E}_x[u_r(X_{t\wedge\tau_{r}})]-u_r(x)=\mathbb{E}_x\left[\int_{0}^{t\wedge\tau_{r}}Lu_r(X_s)\d s\right]\leqslant -\frac{1}{2}\mathbb{E}_x[ t\wedge\tau_{r}],
\end{equation}
where $\tau_r:=\inf\{t\geqslant 0: |X_t|\leqslant r\}.$
By letting $t\rightarrow\infty,$ we obtain 
$$M_r:=\sup_{|x|>r}\mathbb{E}_x[\tau_{r}]\leqslant 2\sup_{|x|>r} u_r(x)=2\left(\frac{1}{r}+\int_{r}^{\infty}\frac{1}{s\tilde{g}(s)} \d s\right)<\infty,$$
hence $\lim\limits_{r\rightarrow\infty}M_r=0$.  
If $1<\alpha<2$, then the dimensional free Harnack inequality holds (see \cite[Corollary 2.2(3)]{ww14}), so  by \cite[Corollary 3.1(2)]{wy11}, the   density $p_t(x,y)$ with respect to $\pi$ exists.   By Lemma \ref{jump heat ker} below, $\|p_s(x,\cdot)\|_{L^2(\pi)}$ is locally bounded, thus $\kappa\geqslant\lambda_1$ by Corollary \ref{heat-ker}. 

Since by \eqref{g}, $\lim_{r\rightarrow\infty}g(r)=\infty$, according to \cite[Theorem 1.1(a)]{wj13}, the process is exponentially ergodic. Thus $\lambda_1>0$. Then we have $\kappa>0$, i.e. the process is uniformly ergodic.
\deprf

\begin{lemma}\label{jump heat ker}
Assume that $1<\alpha<2$ and $p_t(x,y)$ is the transition density, then for any $s>0$ and $x\in\mathbb{R}$, there exists a constant $C>0$ such that for any $ r>0,$
\begin{equation}
\|p_s(x,\cdot)\|_{L^2(\pi)}^2\leqslant \frac{1}{\pi(B(x,r))}\e^{V_s(r)}, 
\end{equation}
where  $$V_s(r)=\frac{ 2Cr^{2}}{(s \wedge 1)^{\frac{2}{\alpha}}}+\frac{C(2r^2)^{\frac{\alpha}{2(\alpha-1)}}}{(s \wedge 1)^{\frac{1}{ \alpha-1}}}.$$ 
\end{lemma}
\prf
According to Harnack inequality (see \cite[Theorem 2.1]{ww14}), 
for any $s>0, \ x, y \in \mathbb{R}^{d}$ and positive $f \in \mathscr{B}_{b}\left(\mathbb{R}^{d}\right)$,
$$
\left(P_{s} f(y)\right)^{2} \leqslant\left(P_{s} f^{2}(x)\right)\e^{V_s(|x-y|)}.
$$
Now by choosing $f(y)=n\wedge p_s(x,y)$, the desired result follows from  a similar argument to the proof of Lemma \ref{non-sym density}.

\deprf


\section{Estimate of $\kappa$ by hitting time}\label{esti by hit}
Now we are going to another direction, for seeking the lower bound of $\kappa$ by the uniform moment of hitting time to some bound set. In this case, we will first obtain the lower bound of $\lambda$ (or $\lambda_1$) by using the hitting time.

This strategy was done well for the reversible Markov chain on countable state space, see for example \cite{myh10}. 

Let $X_t$ be a continuous-time Markov chain on a denumerable state space $E$. The transition function $P_t=(p_{ij}(t))_{i,j\in E}$ is reversible with respect to the stationary distribution $\pi=(\pi_i)_{i\in E}$:
$$
\pi_ip_{ij}(t)=\pi_jp_{ji}(t),\q i,j\in E, t\geqslant0.
$$
Let $\tau_x=\inf\set{t\geqslant0:X_t=x}$ be the hitting time to state $x\in E$. The following lemma can be found in \cite[Proposition 3.2]{cmf00'}.

\bg{lemma}\lb{lmd0}
Let $P^x(t)=(p_{ij}^x(t))_{i,j\not=x}$ be the killed semigroup upon $x\in E$:
$$
p_{ij}^x(t)=\mathbb P_i[X_t=j,t<\tau_x].
$$
Then $\lmd_1\geqslant\lmd^x$, where  $\lmd^x$ is the Dirichlet eigenvalue of killed process upon $x$:
$$
\lmd^x=-\lim_{t\rar\ift}\fr1t\log\|P_t^x\|_{2\rightarrow 2}
$$
\end{lemma}

\begin{theorem}\lb{mc1}
Under the above assumptions, it holds that
$$
\kp\geqslant\sup_{x\in E}\(\sup_{i\not=x}\mathbb E_i\tau_x\)^{-1}.
$$
\end{theorem}

\prf To apply Corollary \ref{heat-ker}, by Lemma \ref{lmd0}  we  only need to prove that
$$
\lmd^x\geqslant \(\sup_{i\not=x}\mathbb E_i\tau_x\)^{-1}.
$$
Assume  $ M_{x}:=\sup_{i\in E}\mathbb{E}_i\tau_{x}<\infty$.
By \eqref{exp mom}, for any $\beta<1/ M_{x}$, 
$$
\sup_{i\in E}\mathbb{E}_i e^{\bt \tau_{x}}\leqslant \fr1{1-\bt M_x}.
$$
So
\begin{equation}
\|P_t^x\|_{\infty\rightarrow \infty}=\sup_{i\in E}\sum_{j\in E}p^x_{ij}(t)=\sup_{i\in E}\mathbb{P}_i[t<\tau_{x}]\leqslant\frac{1}{1-\beta M_0}\e^{-\beta t}.
\end{equation}
By the symmetry $\|P_t^x\|_{1\rightarrow 1}=\|P_t^x\|_{\infty\rightarrow \infty}$, the interpolation theorem (cf. \cite{RS78}) implies 
\begin{equation}
\|P_t^x\|_{2\rightarrow2}\leqslant\frac{1}{1-\beta M_x}\e^{-\beta t}\q\text{for}\ 0<\beta<1/ M_{x}.
\end{equation}
Hence 
$$\lambda^x\geqslant\frac{1}{ M_{x}}.
$$
Consequently, we have $\kp\geqslant \sup_{x\in E}1/M_x$.
\deprf


We can also use this strategy to one-dimensional reversible Markov processes. 
Specially, as an example, we consider  time-changed symmetric $\alpha$-stable processes.

Let $X$ be a symmetric $\alpha$-stable processes on $\mathbb{R}$ with generator $-(-\Delta)^{\alpha/2}$, $\alpha\in(1,2)$,
where $-(-\Delta)^{\alpha/2}$ is the fractional Laplacian operator.
Note that this process is recurrent but not ergodic. 

Let $a$ be a positive function  so that $1/a$ is $L^1(\mathbb{R}; dx) $ locally integrable. 
Consider the following process $Y=(Y_t)_{t\geqslant 0}$:
$$Y_t:=X_{T_t}, \q \text{where}\q T_t=\inf\left\{s\geqslant0: \int_{0}^{s}a(X_u)^{-1}\d u>t\right\}.
$$
We say $Y$ is a time-changed $\alpha$-stable process, which remains recurrent (cf. \cite[Theorem 5.2.5]{CM12}). By \cite{CW14}, $Y$ is a symmetric strong Markov process with the reversible measure $\pi(\d x)=a(x)^{-1}\d x$, and the associated  regular Dirichlet form $({D},\mathscr{D})$ is given by
\begin{equation}\label{Diri form}
	D(f, g)=\frac{1}{2} \int_{\mathbb{R}} \int_{\mathbb{R}}(f(x)-f(y))(g(x)-g(y)) \frac{C_{\alpha}\mathrm{d} x \mathrm{d} y}{|x-y|^{1+\alpha}} , \ \ f,g\in \mathscr{F},
\end{equation}
where 
\begin{equation}\label{Diri domain}
	\mathscr{D}:=\{u\in L^2(\pi): \ D(u,u)<\infty \}, 
\end{equation}
and  $C_\alpha=\frac{\alpha2^{\alpha-1}\Gamma((\alpha+1)/2)}{\sqrt{\pi}\Gamma(1-{\alpha}/{2})}$.
Since the recurrence implies that $Y$ is Lebesgue irreducible (see \cite[Page 42]{chung86} for the definition),
by \cite[Theorem 4.1.1 and Theorem 4.2.1]{chung86}, if $\pi(\mathbb{R}) <\infty$, then $Y$ is ergodic, i.e. for any  $x\in\mathbb{R}$,
$\lim_{t\rightarrow\infty}\|P_t(x,\cdot)-\pi\|_{\mathrm{Var}}=0.$


\begin{theorem}\label{scstable}
For the time-changed process $Y$, if $I:=\int_{\mathbb{R}}a(x)^{-1}|x|^{\alpha-1}\d x<\infty$,
then  the process is uniformly ergodic and
$$
\kappa\geqslant\frac{1}{\omega_{\alpha}I}>0,
$$
where 
$$\omega_{\alpha}:=-\frac{1}{\cos(\pi\alpha/2)\Gamma(\alpha)}>0.
$$
\end{theorem}
\prf  Let $\tau_0=\inf\{t\geqslant0:X_t=0\}$. Define the killed transition semigroup $P_t^0$ by $$P_t^0(x,A)=\mathbb{P}_x[X_t\in A, t<\tau_0]\q  \ \text{for}\q A\in \mathscr{B}(\mathbb{R}),
$$
and  the Green function $G_0^X(\cdot,\cdot)$ for $X$ killed upon $0$ by
$$
G_0^X(x,\d y)=\int_0^{\infty}P_t^0(x,\d y)\d t=G_0^X(x,y)\d y,
$$
where  $G_0^X(x,y)$ is the Green function for $X$ killed upon $0$ (cf. \cite[Page 152]{Get20}):
$$G_0^X(x,y)=-\frac{1}{2\Gamma(\alpha)\cos\left(\frac{\pi\alpha}{2}\right)}\left(|y|^{\alpha-1}+|x|^{\alpha-1}-|y-x|^{\alpha-1}\right).$$
By \cite[(4.25)]{DK20},  we can represent the Green function $G_0^Y(\cdot,\cdot)$ for $Y$ killed upon $0$ as
$$
G_0^Y(x,A)=\int_A G_0^X(x,y)a(y)^{-1}\d y.
$$
Therefore, 
$$\mathbb{E}_x\tau_0^Y=\int_0^\ift\mathbb E_x\mathbf{1}_{\mathbb R}(Y_t^0)\d t=\int_{\mathbb{R}}G_0^Y(x,\d y)=\int_{\mathbb{R}}G_0^X(x,y)  a(y)^{-1}\d y,
$$
where $\tau_0^Y=\inf\{t\geqslant0:Y_t=0\}$.
According to Lemma \ref{trian} below,  we have
$$
M_0^Y:=\sup_x\mathbb{E}_x\tau^Y_0\leqslant -\frac{1}{\Gamma(\alpha)\cos\left(\frac{\pi\alpha}{2}\right)}\int_{\mathbb{R}}|y|^{\alpha-1}  a(y)^{-1}\d y=\omega_{\alpha}I^{\sigma,\alpha},
$$
By \cite[Lemma 3.2]{GT12}, 
$\lambda_0\geqslant \(M_0^Y\)^{-1},$
where $$\lambda_0:=\inf\{D(f,f):f\in\mathscr{D},\pi(f^2)=1,f(0)=0\},$$
and $(D,\mathscr{D})$ is the Dirichlet form of $Y$ given by \eqref{Diri form} and \eqref{Diri domain}. 
It is well known that 
$\lambda_1\geqslant\lambda_0 \ ( \text{see \cite[Proposition 3.2]{cmf00'}})$, 
thus $\lambda_1\geqslant (M_0^Y)^{-1}.$ Now our result follows by Corollary \ref{heat-ker}.
\deprf

\begin{lemma}\label{trian}
For any $x,y\in\mathbb{R}$ and $\alp\in(1,2)$, 
$$|y|^{\alpha-1}+|x|^{\alpha-1}-|y-x|^{\alpha-1}\leqslant2(|x|\wedge|y|)^{\alpha-1},
$$
\end{lemma}
\prf
Let $a=|x\wedge y|,$ $b=|x\vee y|$. Then $|x|\wedge|y|=a\wedge b$.

(1) When $xy=0$, it is trivial. 

(2) When $xy<0$,  
$$|y|^{\alpha-1}+|x|^{\alpha-1}-|y-x|^{\alpha-1}=a^{\alpha-1}+b^{\alpha-1}-(b+a)^{\alpha-1}\leqslant  (a\wedge b)^{\alpha-1}=(|x|\wedge|y|)^{\alpha-1}.
$$

(3) When $xy>0$, we only need to consider $x,y>0$. Note that for any $c_1,c_2>0$,
$$
(c_1+c_2)^{\alpha-1}\leqslant c_1^{\alpha-1}+c_2^{\alpha-1}.
$$
Therefore, $b^{\alpha-1}\leqslant (b-a)^{\alpha-1}+a^{\alpha-1}$, which means that $a^{\alpha-1}+b^{\alpha-1}-(b-a)^{\alpha-1}\leqslant2a^{\alpha-1}.$ 
Then
$$
y^{\alpha-1}+x^{\alpha-1}-|y-x|^{\alpha-1}\leqslant2(x\wedge y)^{\alpha-1}.
$$
\deprf
\begin{remark} For one-dimensional  time-changed symmetric $\alpha$-stable process with $\alpha\in(1,2)$, \cite[Theorem 1.7]{CW14} proves the following sufficient condition for uniform ergodicity:
$$
\liminf_{|x|\rightarrow\infty}\frac{a(x)^{1/\alpha}}{|x|^\gamma}>0\ \text{for some} \ \gamma>1.
$$
Theorem \ref{scstable} is an extension of this result.
\end{remark}

{\bf Acknowledgement}\
This work is supported in part by the  National Key Research and Development Program of China (2020YFA0712900), the National Natural Science Foundation of China (Grant No.11771047 and No.11771046) and the project from the Ministry of Education in China.


\bibliographystyle{plain}
\bibliography{convergence_rate}

\begin{thebibliography}{10}

\bibitem{Ander91}
W.~Anderson.
\newblock {\em {Continuous-time Markov chains}}.
\newblock Springer-Verlag, New York, 1991.

\bibitem{Bax05}
P.H. Baxendale.
\newblock {Renewal theory and computable convergence rates for geometrically
  ergodic Markov chains}.
\newblock {\em Ann. Appl. Probab.}, 15:700--738, 2005.

\bibitem{brw01}
V.I. Bogachev, M.~Rockner, and F.Y. Wang.
\newblock Elliptic equations for invariant measures on finite and infinite
  dimensional manifolds.
\newblock {\em J. Math. Pure. Appl.}, 80(2):177--221, 2001.

\bibitem{cmf99}
M.-F. Chen.
\newblock {Analytic proof of dual variational formula for the first eigenvalues
  in dimension one.}
\newblock {\em Sci. in China Ser. A}, 42:805–--815, 1999.

\bibitem{cmf00'}
M.-F. Chen.
\newblock {Explicit bounds of the first eigenvalues.}
\newblock {\em Sci. in China Ser. A}, 43:1051--1059, 2000.

\bibitem{cmf04}
M.-F. Chen.
\newblock {\em {From {M}arkov chains to non-equilibrium particle systems, 2nd
  edition}}.
\newblock World Scientific, Singapore, 2004.

\bibitem{cmf05}
M.-F. Chen.
\newblock {\em {Eigenvalues, inequalities, and ergodic theory}}.
\newblock Springer, London, 2005.

\bibitem{cmf95}
M.-F. Chen and F.-Y. Wang.
\newblock Estimation of the first eigenvalue of the second order elliptic
  operators.
\newblock {\em J. Funct. Anal.}, 131(2):345--363, 1995.

\bibitem{cmf97}
M.-F. Chen and F.-Y. Wang.
\newblock Estiamtion of spectral gap for elliptic operators.
\newblock {\em Trans. Amer. Math. Soc.}, 349(3):1239--1267, 1997.

\bibitem{CM12}
Z.-Q. Chen and M.~Fukushima.
\newblock {\em {Symmetric Markov processes, time change, and boundary theory}}.
\newblock Princeton Univ, Princeton. Press, 2012.

\bibitem{CW14}
Z.-Q. Chen and J.~Wang.
\newblock {Ergodicity for time-changed symmetric stable processes}.
\newblock {\em Stoch. Proc. Appl.}, 124(9):2799--2823, 2014.

\bibitem{chung86}
K.L. Chung.
\newblock Doubly-feller process with multiplicative functional.
\newblock {\em In Seminar on Stochastic Processes, Birkh\"{a}user.}, 12:63--78,
  1986.

\bibitem{cz95}
K.L. Chung and Z.~Zhao.
\newblock {\em {From Brownian motion to Schr\"{o}dinger's equation}}.
\newblock Springer, New York, 1995.

\bibitem{DK20}
L.~Doering and A.E. Kyprianou.
\newblock Entrance and exit at infinity for stable jump diffusions.
\newblock {\em Ann. Probab.}, 48(3):1220--1265, 2020.

\bibitem{DMPS18}
R.~Douc, E.~Moulines, P.~Priouret, and P.~Soulier.
\newblock {\em {Markov Chains}}.
\newblock Springer Series in Operations Research and Financial Engineering.
  Springer, Cham, 2018.

\bibitem{MT95}
D.~Down, S.P. Meyn, and R.T. Tweedie.
\newblock {Exponential and uniform ergodicity of Markov processes}.
\newblock {\em Ann. Probab.}, 23:1671--1691, 1995.

\bibitem{FLZ20}
C.~Foucart, P.-S. Li, and X.~Zhou.
\newblock {On the entrance at infinity of Feller processes with no negative
  jumps}.
\newblock {\em Stat. Probab. Lett.}, 165, 2020.

\bibitem{Get20}
R.~K. Getoor.
\newblock Continuous additive functionals of a markov process with applications
  to processes with independent increments.
\newblock {\em J. Math. Anal. Appl.}, 13:132--153, 1966.

\bibitem{GT12}
A.~Grigor'yan and A.~Telcs.
\newblock Two-sided estimates of heat kernels on metric measure spaces.
\newblock {\em Ann. Probab.}, 40:1212--1284., 2012.

\bibitem{IM65}
K.~It\^{o} and H.P.~Mckean Jr.
\newblock {\em {Diffusion Processes and Their Sample Paths}}.
\newblock Springer, Berlin, 1965.

\bibitem{myh02}
Y.-H. Mao.
\newblock {Strong ergodicity for Markov processes by coupling methods}.
\newblock {\em J. Appl. Prob.}, 39:839--852, 2002.

\bibitem{myh06}
Y.-H. Mao.
\newblock {Convergence rates in strong ergodicity for Markov processes}.
\newblock {\em Stoch. Proc. Appl.}, 116:1964--1976, 2006.

\bibitem{myh10}
Y.-H. Mao.
\newblock {Convergence rates for reversible Markov chains without the
  assumption of nonnegative definite matrices}.
\newblock {\em Sci. China Math.}, 53(8):1979–--1988, 2010.

\bibitem{mz17}
Y.-H. Mao and C.~Zhang.
\newblock { Uniform convergence rates for birth and death processes}.
\newblock {\em Markov Proc. Relat. Fields.}, 23:467--483, 2017.

\bibitem{MT09}
S.P. Meyn and R.L. Tweedie.
\newblock {\em {Markov chains and stochastic stability, 2nd edition}}.
\newblock {Cambridge University Press, New York}, 2009.

\bibitem{RS78}
M.~Reed and B.~Simon.
\newblock {\em {Methods of Modern Mathematical Physics, Vol. II,}}.
\newblock Academic Press, New York, 1978.

\bibitem{rosen95}
J.S. Rosenthal.
\newblock {Minorization conditions and convergence rates for Markov chain Monte
  Carlo }.
\newblock {\em J. Amer. Statist. Assoc.}, 90:558–--566, 1995.

\bibitem{sen88}
E.~Seneta.
\newblock {Perturbation of the stationary distribution measured by ergodicity
  coefficients }.
\newblock {\em Adv. Appl. Probab.}, 20:228--230, 1988.

\bibitem{ST89}
A.D. Sokal and L.E.Thomas.
\newblock {Exponential convergence to equilibrium for a class of random-walk
  models}.
\newblock {\em J. Stat. Phys.}, 54:797–--828, 1989.

\bibitem{wfy05}
F.-Y. Wang.
\newblock {\em {Functional inequalities, Markov semigroups, and spectral
  theory}}.
\newblock Science press, Beijing, 2005.

\bibitem{wfy13}
F.-Y. Wang.
\newblock {\em {Analysis for diffusion processes on Riemannian manifolds.}}
\newblock World scientific, Singapore, 2013.

\bibitem{ww14}
F.-Y. Wang and J.~Wang.
\newblock {Harnack inequality for stochastic differential equations driven with
  L\'{e}vy noises,}.
\newblock {\em J. Math. Anal. Appl.}, 410:513--523, 2014.

\bibitem{wy11}
F.-Y. Wang and C.~Yuan.
\newblock {Harnack inequalities for functional SDEs with multiplicative noise
  and applications}.
\newblock {\em Stochastic Process. Appl.}, 121:2692–--2710., 2011.

\bibitem{WJ09}
J.~Wang.
\newblock {Logarithmic Sobolev inequality and strong ergodicity for birth-death
  processes.}
\newblock {\em Front. Math. China}, 4:721--726, 2009.

\bibitem{wj13}
J.~Wang.
\newblock {Exponential ergodicity and strong ergodicity for SDEs driven by
  symmetric $\alpha$-stable processes.}
\newblock {\em Appl. Math. Lett.}, 26:654--658, 2013.

\bibitem{zyh18}
Y.-H. Zhang.
\newblock {Criteria on ergodicity and strong ergodicity of single death
  processes.}
\newblock {\em Front. Math. China}, 13(5):1215–--1243, 2018.

\end{thebibliography}

\end{document}